\documentclass[10pt,a4paper]{amsart}
\usepackage{amsfonts}
\usepackage{mathtools}
\usepackage{amsthm}
\usepackage{amsmath}
\usepackage{amssymb}
\usepackage{amscd}
\usepackage[ngerman,english]{babel}
\usepackage[mathscr]{eucal}
\usepackage{indentfirst}
\usepackage{graphicx}
\usepackage{graphics}
\usepackage{comment}
\usepackage{pict2e}
\usepackage{epic}
\numberwithin{equation}{section}
\usepackage[margin=3.5cm]{geometry}
\usepackage{epstopdf} 
\usepackage{color}
\usepackage[hidelinks]{hyperref}
\usepackage[numbers]{natbib}

\theoremstyle{plain}
\newtheorem{Th}{Theorem}[section]

\newtheorem{Cor}[Th]{Corollary}

\DeclareMathOperator{\R}{\mathbb{R}}

\DeclareMathOperator{\N}{\mathbb{N}}
\DeclareMathOperator{\di}{\text{div}}
\DeclareMathOperator{\supp}{\text{supp}}
\DeclareMathOperator{\tr}{\text{tr}}
\DeclareMathOperator{\codim}{\text{codim}}
\DeclareMathOperator{\dist}{\text{dist}}

\theoremstyle{definition}

\newtheorem{Rem}[Th]{Remark}
\newtheorem{?}[Th]{Problem}

\newcommand{\norm}[2]{\left\lVert #1 \right\rVert_{#2}}
\newcommand{\f}[2]{\frac{#1}{#2}}

\begin{document}

\title[Energy bounds for an equation related to wave maps ]{Energy bounds  for a fourth-order equation  in low dimensions related to wave maps \\  }

\author[T.Schmid]{Tobias Schmid}

	\address{EPFL SB MATH PDE,
	B\^atiment MA,
	Station 8,
	CH-1015 Lausanne, Switzerland}
\email{tobias.schmid@epfl.ch}



\thanks{The author acknowledges funding by the Deutsche Forschungsgemeinschaft (DFG, German Research Foundation) – Project-ID 258734477 – SFB 1173}

\begin{abstract} For compact, isometrically embedded Riemannian manifolds $ N \hookrightarrow \R^L$, we introduce a  fourth-order version of the  wave maps equation. 
By energy estimates, we prove an \emph{a priori} estimate for smooth local solutions in the energy subcritical dimension $ n = 1,2$. The estimate excludes blow-up of a Sobolev norm in finite existence times. In particular, combining this with recent work of local well-posedness of the Cauchy problem, it follows that for smooth initial data with compact support, there exists a (smooth) unique global solution in dimension $n =1,2$. We also give a proof of the uniqueness of solutions that are bounded in these Sobolev norms.
\end{abstract}

\maketitle

\section{Introduction}
Let $ (N,h)$ be a (compact) Riemannian manifold, isometrically embedded (by Nash's theorem) into euclidean space $N \hookrightarrow \R^L$. For a Riemannian manifold $(M,g) $, we introduce the action functional
$$ \mathcal{L}(u) = \f{1}{2}\int_{0}^T \int_M |\partial_tu(x,t)|^2 - | \Delta_{g(x)} u (x,t)|^2~ dV_g(x)dt,~~~~dV_g = \sqrt{\det g}~dx$$
for (smooth) maps $u : M \times [0,T) \to N $. The action $ \mathcal{L}$ is critical with respect to $u$ if 
$$ \f{d}{d \delta} \mathcal{L}(u + \delta \Phi)_{|_{\delta = 0}} = 0,~~ \Phi \in C^{\infty}_c( M \times (0, T), \R^L),~ \Phi(x,t) \in T_{u(x,t)}N,~~ (x,t) \in M \times [0,T).$$
In this case, $ u$ satisfies the condition
\begin{equation}\label{basic}
\partial_{t}^2u(x,t) + \Delta_g^2u(x,t) \perp T_{u(x,t)}N, ~ (x,t) \in M \times (0,T),
\end{equation}
where $ \Delta_g $ denotes the Laplace-Beltrami operator on $(M,g)$. 
 More explicitly, we use the fact that there exists a smooth familiy of orthogonal (linear) projector
$$ P_p : \R^L \to T_pN ,~p \in N, $$
in order to expand \eqref{basic} into the Euler-Lagrange equation
\begin{align}\label{EL}
\partial_t^2 u + \Delta_g^2 u =& (I - P_u)( \partial_t^2 u + \Delta_g^2 u)\\ \nonumber
=& ~dP_u(\partial_t u, \partial_t u) +  \Delta_g ( \tr_g dP_u( \nabla u, \nabla u))\\[3pt] \nonumber  
&+ 2\di_g( dP_u(\nabla u, \Delta_g u)) - dP_u( \Delta_g u , \Delta_g u).
\end{align}
The right hand side of \eqref{EL} is highly non-linear and we consider \eqref{EL} to be a prototypical model of a geometric fourth-order equation. In particular, we expect the long-time behaviour, such as scattering v.s. blow up, to be influenced by the choice of $M, N $ similarly as e.g. for the well-known Schr\"odinger maps or wave maps equation. For the latter we refer to \citep{shatahstruwe} for an overview.
 The operator 
$ \partial_t^2 + \Delta^2 $
typically relates to the Euler-Bernoulli beam equation or two-dimensional elastic plate equations (under idealized assumptions). In general the potential energy, say in case $M = \R^n $,
$$ \mathcal{E}(u) = \f12 \int_{\R^n} |\Delta u |^2~dx,~~ \Delta u = \partial_i \partial^i u, $$
approximates the bending energy of an elastic, stiff membrane with small deflections.  
In the following, we give  an existence argument based on global energy bounds for \eqref{EL} in the case $ M = \R^n$ in \emph{low} dimensions $ n \in \{1,2\}$.\\[3pt]
The equation \eqref{EL} is Hamiltonian with energy 
\begin{equation}\label{cons}
E(u(t)) = \f{1}{2} \int_{\R^n} | u_t|^2 + | \Delta u|^2 dx.
\end{equation}
Further \eqref{EL} has scaling $ u_{\lambda}(t,x) = u(\lambda^2 t, \lambda x)$, i.e.  if $u$ solves \eqref{EL} on $[0,T) \times \R^n$,  
then $ u_{\lambda }$ for $ \lambda > 0$ is also a solution on $ [0,\lambda^{-2}T) \times \R^n$. 
Since for local solutions
\[
E(u_{\lambda}(t) ) = \lambda^{4 -n} E(u(\lambda^2t)),~~ \lambda > 0,
\]
we refer to dimensions $ n < 4$ as energy subcritical. Here one heuristically expects the conservation of \eqref{cons} to be accessible for  bounding large data solutions of \eqref{EL} (with sufficient regularity).
\begin{Rem}
	The projector maps are  derivatives of the metric distance (with respect to $N$) in $\R^L$, ie.
	$$ p = \Pi(p) + \f{1}{2} \nabla_p ( \dist^2(p, N)),~~ P_p = \nabla_p \Pi(p),~~~\dist(p,N) < \delta_0.$$ 
	We note that via this representation, it is possible to extend this family smoothly to all of $ \R^L$ in order to solve the Cauchy problem for \eqref{EL} without restricting the coefficients a priori.
\end{Rem}
~~\\
We now consider the Cauchy problem on $ M = \R^n $
\begin{align}\label{CP}
\left\{
\begin{matrix*}[l]
\partial_{t}^2u(t,x) + \Delta^2u(t,x) \perp T_{u(t,x)}N, & (t,x) \in (0,T) \times \R^n\\[6pt]
(u(0,x),u_t(0,x)) = ( u_0(x), u_1(x)),  & x \in \R^n \\[6pt]
(u_0,u_1): \R^n \to TN,~~ u_1(x) \in T_{u_0(x)}N & x \in \R^n
\end{matrix*}
\right.
\end{align}
For \eqref{CP}, we obtain the following result.
\begin{Th}\label{main}
	Let $ n \in \{  1,2\}$ and $ u \in C^{\infty}(\R^n \times [0, T), N)$ be a local solution of \eqref{CP}. Assume further 
	$$ u - u_0 \in C^0([0, T), H^{n+2}(\R^n))\cap C^1([0, T), H^n(\R^n)).$$
	Then there holds 
	\begin{equation}
	\limsup_{t \nearrow T}( \norm{ u_t(t)}{H^{n}} + \norm{\nabla u(t)}{H^{n+1}} ) < \infty,
	\end{equation}
	as long as $ T < \infty $.
\end{Th}
Recently in  \citep{LammSchnaubeltHerrSchmid}, the authors proved local wellposedness (in high regularity) and a blow up condition for the Cauchy problem \eqref{CP}, which (by the proof of Theorem \ref{main}) implies the following 
\begin{Cor}\label{Cor}
	Let $ n \in \{  1,2\}$ and $ u_0, u_1 : \R^n \to \R^L,~ u_0(x) \in N,~ u_1(x) \in T_{u_0(x)}N$ for $ x \in \R^n $ and such that 
$$ (\nabla u_0, u_1) \in H^{k}(\R^n)\times H^{k-1}(\R^n),$$
for $ k \in \N $ with $ k \geq  n+1$.
	Then the Cauchy problem \eqref{CP} has a global solution  $ u : \R^n \times \R   \to N $ with
	$$ u - u_0 \in C^0(\R, H^{k+1}(\R^n))\cap C^1(\R, H^{k-1}(\R^n)).$$
	In particular, if $ u_0, u_1 $ are smooth and  $ \supp(\nabla u_0), \supp(u_1) $ are compact, then there exists a global smooth solution of \eqref{CP}. 
\end{Cor}
~~\\
The work is part of the authors PhD thesis \citep{Sschmid} and we conclude this introduction with a few remarks.
In the sense explained above, \eqref{basic} and \eqref{EL} are higher order versions of the wave maps equation
		\begin{equation} \label{wavemap}
		\square_g u = dP_u(\partial_tu , \partial_tu) - \tr_g dP_u( \nabla u, \nabla u),
		\end{equation}
		with the d'Alembert operator $ \square_g = \partial_t^2 - \Delta_g $.
		Equation \eqref{wavemap} is the Euler Lagrange equation of the action functional
		 $$ \mathcal{L}(u) = \int_0^T \int_M L(u)~dV_g ~dt$$
		  on the Riemannian manifold $ (M,g)$ with Lagrangian $ L(u) = \f{1}{2} \tilde{g}^{\alpha \beta}\langle \f{\partial u}{\partial x_{\alpha}},   \f{\partial u}{\partial x_{\beta}}\rangle$, and where  $ \f{\partial}{\partial x_0} = \f{\partial}{\partial t}$ and $ \tilde{g} = - dt^2 + g$. This wave equation has been studied intensively in the past, especially as a model problem for nonlinear dispersion and singularity formation. We refer to \citep{shatahstruwe} and \citep{grillakisgeba} for an overview over the wellposedness and singularity theory of the Cauchy problem for the wave maps equation \eqref{wavemap}. For \eqref{wavemap}, the action functional $ \mathcal{L}$ is independent of the embedding $ N \hookrightarrow \R^L $. In our case however, there is an intrinsic map equation, arising from critical points of the (embedding independent) functional
		$$ \mathcal{L}_i(u) =  \int_0^T\int_{M} | \partial_t u|_h^2 - | \tr_g(\nabla du)  |_h^2~dV_g~dt. $$
		where $ \nabla $ denotes the Levi-Civita connection of the pullback bundle $ u^* TN $ endowed with the pullback metric $ u^*h$ and the energy potential is given by the tension field $ \tau_g(u) = \tr_g (\nabla du)$ of $u$. Moreover, first variations are calculated intrinsically as follows.
		$$ \frac{d}{d \delta} \mathcal{L}_i (u^{\delta})_{|_{\delta =0}} = 0,$$
		where $ u^{\delta}  \in C^{\infty}( M \times [0,T),N)$ is a family of  maps depending smoothly on $ \delta \in (- \delta_0, \delta_0)$ with $ u = u^0$  and 
		compact $ \supp(u - u^{\delta}) \subset M \times (0,T)$ for $ |\delta | < \delta_0$. \\[3pt]
Then the Euler-Lagrange equation, which has been calculated for static solutions e.g. in \citep{jiang}, becomes
\begin{equation} \label{EL2}
\nabla_t \partial_t u + \Delta^2_{g,h} u + R(u)( d u, \Delta_{g,h} u) d u = 0,
\end{equation}
		where $R$ is the curvature tensor and in the covariant notation, we set $ \Delta_{g,h} u = \tr_g(\nabla du)$, and use $\Delta^2_{g,h} u = \Delta^{}_{g,h}(\Delta_{g,h} u) = \tr_g(\nabla \tau_g(u))$.\\[10pt] 
Static solutions of \eqref{EL} (and \eqref{EL2}) are \emph{extrinsic (and intrinsic) biharmonic maps}, i.e. they are maps $ u : (M,g) \to (N,h) $ between Riemannian manifolds that are critical for the (intrinsic or extrinsic) energy functional
$$ F(u) = \frac{1}{2}\int_M | \tr_g(\nabla du) |^2_h~ dV_g,~~~~~E(u) = \frac{1}{2}\int_M | \Delta_g u |^2_h ~dV_g,~\text{respectively}$$
where the latter is defined subject to an isometric embedding $ (N,h) \hookrightarrow \R^m.$ Biharmonic maps (resp. the Euler Lagrange equation of $ E$ and $F$) and their heat flows have been studied intensively in the past.\\[10pt]
\section{Related work and local wellposedness in high regularity}~~\\
In \citep{LammSchnaubeltHerr}, the authors pove the existence of a global weak solution into round spheres $ \mathbb{S}^{L-1} \subset \R^L$. This is done by a penalization functional of Ginzburg Landau type, which then gives a uniform energy bound in the penalty parameter. To prove convergence of such approximations, the authors depend on the geometry of the sphere, more precisely, the equation can be rewritten in divergence form. This argument has been used for the wave maps equation \eqref{wavemap} with $ N = \mathbb{S}^{L-1}$ and $ M = \R^n $  in \citep{shatah1988weak} and further the divergence form has been used in \citep{Strzelecki}, in order to prove weak compactness of the class of stationary solutions of \eqref{EL} on the domain $M = \R^4$.\\[4pt]
 In \citep{schmid2}, the author proved the existence of low-regularity and global smooth solutions of \eqref{EL}  with small initial data in a scaling critical (Besov) space on $ \R^{n},~ n \geq 3$. The main aspect of \eqref{EL} used in  \citep{schmid2}, is the exploitation of a non-resonant form of the nonlinearity, which essentially resembles the geometric condition \eqref{CP}. Spaces of initial data $D$ are called (scaling) critical, if $ \| (u(0), \partial_t u(0))\|_D = \| (u_{\lambda}(0), \partial_tu_{\lambda}(0))\|_D$ where $\lambda > 0$. For data of regularity close to that assumed in such spaces, the wellposedness of \eqref{EL} is heuristically expected to be at a threshold and thus  particularly challenging.\\[4pt]
As mentioned above, in the recent work \citep{LammSchnaubeltHerrSchmid}, the authors prove local wellposedness of the Cauchy problem \eqref{CP}. More precisely, let $u_0, u_1 : \R^n \to \R^L,~ u_0(x) \in N,~ u_1(x) \in T_{u_0(x)}N$ for $\mathcal{L}^n$ a.e. $ x \in \R^n $ with
$$ (\nabla u_0, u_1) \in H^{k-1}(\R^n)\times H^{k-2}(\R^n),~~ k > \left \lfloor \f{n}{2} \right \rfloor +2,~ k \in \N.$$
Then there exists a $ T > 0 $ and a (unique) solution  $ u : \R^n \times [0,T) \to N $ of \eqref{CP} with
$$ u-u_0 \in C^0([0,T), H^k(\R^n)) \cap C^1([0,T), H^{k-2}(\R^n)).$$
From this, we note that in particular we obtain Corollary \ref{Cor} from a blow up condition contained in \citep{LammSchnaubeltHerrSchmid}. In the following, recall that the energy functional
\begin{equation*}
E(u(t)) = \f{1}{2} \int_{\R^n} | u_t|^2 + | \Delta u|^2 dx,
\end{equation*}
is formally conserved along solutions $u$. This implies the bound
\begin{equation}\label{grad}
\f{d}{dt} \int_{\R^n} | \nabla u|^2dx \leq E(u(0)).
\end{equation}
Both, the conservation of \eqref{cons} and \eqref{grad}, will be used in the following for smooth solutions. We further note that below in section \ref*{uniqueness}, we include a short argument for the uniqueness of such solutions in dimension $ n = 1,2,3$.
\section{Proof of Theorem \ref{main}}
Since for solutions $u$ of \eqref{EL}, resp. the Cauchy problem \eqref{CP}, the term $ \partial_{t}^2u + \Delta^2u $ is a section over the normal bundle of $u^*(TN)$, we let $\codim(N) = L-l $ for $ l \in \N,~ l \leq L $ and first assume the normalbundle $ T^{\perp}N $ of $N \subset \R^L$ is parallelizable. This means there exists a frame of (smooth) orthogonal  vectorfields $ \{ \nu_1(p), \dots, \nu_{L-l}(p)\} \subset \R^L,~ p \in N $ with $ \nu_i (p) \perp T_pN $ for every $ p \in N $.\\[5pt] In this case, for any local solution $ u $,  we have an explicit representation for the nonlinearity in terms of $ \nu_i(u) $.
\begin{equation}\label{eq}
\partial_t^2 u + \Delta^2 u =: \sum_{i = 1}^{L-l}G^i(u) \nu_i(u) =: G^i(u) \nu_i(u),
\end{equation}
where $G_i(u) = \langle \partial_t^2 u + \Delta^2 u , \nu_i(u) \rangle $. We thus calculate
\begin{align*}
&\langle \partial_{t}^2 u, \nu_i(u)\rangle = - \langle u_t , d\nu_i(u) u_t \rangle,\\
&\langle \Delta^2 u, \nu_i(u) \rangle = - 3\langle \nabla \Delta u , d\nu_i(u) \nabla u \rangle - \langle \nabla u , d\nu_i(u) \nabla \Delta u \rangle\\
&~~- \langle \nabla u, d^3\nu_i (u)(\nabla u)^3 + 2 d^2 \nu_i(u)( \nabla u, \nabla^2 u ) +  d^2\nu_i(u) (\nabla u, \Delta u) \rangle\\
& ~~ - 2 \langle \nabla^2 u, d^2 \nu_i(u) (\nabla u)^2 + d\nu_i(\nabla^2 u) \rangle - \langle \Delta u , d^2\nu_i(u) (\nabla u )^2 + d \nu_i(\Delta u) \rangle,
\end{align*}
where we denote by  $d^k \nu_i $ the $kth$ order differential of $ \nu_i $ on $N$
 and write $ (\nabla u)^2 ,~ (\nabla u)^3 $ for products of first order derivatives of $u$ with eiter two or three factors, respectively. The precise product, e.g. $ \partial_{x_j} u \cdot \partial^{x_j} u$ or $ \partial_{x_i} u \cdot \partial^{x_j} u \cdot \partial_{x_j} u $, will become clear in the terms of the expansion.
 The result  in Theorem \ref{main} is known for $N = \mathbb{S}^{L-1}$ and $n \leq 2$ thanks to \citep{fan2010regularity}.\\[10pt]
\textit{\underline{Case: n = 2}}
We apply $ \Delta = \partial_i \partial^i $ on both sides of \eqref{eq}. Then, testing the differentiated equation by  $ \Delta u_t $, we infer 
\begin{align}
\f{d}{2dt} \int_{\R^n} (| \Delta u_t|^2 + |\Delta^2 u |^2)dx = \int_{\R^n}\Delta ( G^i(u) \nu_i(u)) \Delta u_t dx. \end{align}
Since $ G^i(u) $ contains derivatives of order three, we can not proceed  by the H\"older inequality. Instead, we follow \citep{fan2010regularity}, where the authors showed that the highest order derivative cancel in the case $ N = \mathbb{S}^{L-1},~ \nu(u) = u $. Since 
\begin{align*}
\Delta ( G^i(u) \nu_i(u)) \Delta u_t = \Delta (G^i(u)) \nu_i(u) \Delta u_t + 2 \nabla (G^i(u)) \cdot\nabla ( \nu_i(u))\Delta u_t + G^i(u) \Delta \nu_i(u) \Delta u_t,
\end{align*}
and
\begin{align*}
0 = \Delta ( \nu_i(u) u_t) = 2  d\nu_i(u)( \nabla u) \cdot \nabla u_t + \nu_i(u) \Delta u_t + d^2 \nu_i(u)(\nabla u)^2  u_t +  d\nu_i(u)(\Delta u) u_t,
\end{align*}
it follows 
\begin{align*}
\Delta ( G^i(u) \nu_i(u)) \Delta u_t =& -\Delta G^i(u) \left( 2 d\nu_i(u) (\nabla u)\cdot \nabla u_t + d^2 \nu_i(u) (\nabla u)^2 u_t  + d\nu_i(u) (\Delta u) u_t\right)\\
& + 2\nabla G^i(u) \cdot  d\nu_i(u) (\nabla u) \Delta u_t\\
& + G^i(u) \left( d^2 \nu_i(u) (\nabla u)^2 + d \nu_i(u) \Delta u \right)\Delta u_t.
\end{align*}
Hence we observe, by integration by parts for the first summand, 
\begin{align*}
\int_{\R^n} \Delta ( G^i(u) \nu_i(u)) \Delta u_t~dx =&~ \int_{\R^n} \nabla G^i(u) \cdot [ 3 d^2 \nu_i(u) (\nabla u)^2 \nabla u_t + 3 d \nu_i(u)(\Delta u) \nabla u_t]dx\\
& + \int_{\R^n} \nabla G^i(u) \cdot [ 4 d\nu_i(u) (\nabla u) \Delta u_t + d^3\nu_i(u) (\nabla u)^3 u_t]dx\\
&  + \int_{\R^n} \nabla G^i(u) \cdot [ 3 d^2 \nu_i(u)(\Delta u, \nabla u) u_t  + d\nu_i(u) (\nabla \Delta u) u_t  ]dx\\
& + \int_{\R^n}  G^i(u) ( d^2 \nu_i(u) (\nabla u)^2 u_t + d \nu_i(u) (\Delta u)) \Delta u_t~dx.
\end{align*}
Instead of deducing bounds for this terms that depend on the normal frame $ \{\nu_1, \dots \nu_{L-l}\}$, we turn to the general case and use the normal projector $ I-P_u : \R^L \to (T_uN)^{\perp} $ along the map $ u: \R^n \times [0,T)  \to N $ in order to represent the nonlinearity in \eqref{eq} as
\begin{equation}
\partial_t^2 u + \Delta^2 u = (I- P_u)(\partial_t^2 u + \Delta^2 u ).
\end{equation}
Here, we proceed similarly, ie. we use 
\begin{align}
\Delta ( (I- P_u)(\partial_t^2 u + \Delta^2 u ) )\Delta u_t = \Delta ( (I- P_u)^2(\partial_t^2 u + \Delta^2 u ) )\Delta u_t,
\end{align}
and hence
\begin{align*}
\Delta ( (I- P_u)^2(\partial_t^2 u + \Delta^2 u ) )\Delta u_t =&  \Delta [( I- P_u)] ((I-P_u)(\partial_t^2 u + \Delta^2 u ) )\Delta u_t\\
& + 2 \nabla (I- P_u) \cdot \nabla ((I-P_u)(\partial_t^2  + \Delta ^2 u )) \Delta u_t\\
& + (\Delta [ (I-P_u)(\partial_t^2 u + \Delta^2 u ) ]) (I - P_u)\Delta u_t.
\end{align*}
In order to treat the last summand, we expand 
$$ 0 = \Delta ( (I - P_u) u_t ) = (I-P_u)\Delta u_t - d^2P_u((\nabla u)^2, u_t) - dP_u(\Delta u, u_t) - 2 dP_u ( \nabla u, \nabla u_t).$$
Hence, as before, integration by parts yields
\begin{align*}
\int_{\R^n} \Delta &( (I-P_u) (\partial_t^2u + \Delta^2 u)) \Delta u_t\\
=& - \int_{\R^n} d^2P_u( (\nabla u)^2, (I-P_u) (\partial_t^2u + \Delta^2 u)) \Delta u_t~dx\\
& - \int_{\R^n} dP_u(\Delta u, (I-P_u) (\partial_t^2u + \Delta^2 u)) \Delta u_t~dx\\
& - 2\int_{\R^n} dP_u(\nabla u, \nabla [ (I-P_u) (\partial_t^2u + \Delta^2 u)]) \Delta u_t~dx\\
& - \int_{\R^n} \nabla [ (I-P_u) (\partial_t^2u + \Delta^2 u)] \cdot \nabla [ dP_u(\Delta u, u_t) + 2 dP_u(\nabla u, \nabla u_t) + d^2P_u((\nabla u)^2, u_t)]~dx.
\end{align*}
We first note the pointwise bounds 
\setlength\jot{.2cm}
\begin{alignat}{2}
&|(I-P_u)(\partial_t^2 u + \Delta^2 u) | \lesssim&&  | u_t|^2 + | \nabla^2 u |^2+  | \nabla^2 u| |\nabla u|^2 + | \nabla^3  u| | \nabla u| + |\nabla u |^4 \label{one} \\ 
&| \nabla [( I- P_u)( \partial_t^2 u + \Delta^2 u )]| \lesssim&& ~~~ |\nabla u_t| | u_t| + | \nabla u | | u_t|^2 + |  \Delta^2 u| | \nabla u |\label{two}\\ \nonumber
&  && + | \nabla^3 u | ( | \nabla^2 u | + | \nabla u|^2 ) + |\nabla u| | \nabla^2 u |^2 + | \nabla u|^3  | \nabla^2 u | + | \nabla u|^5,
\end{alignat}
where the constants only depend on the supremum norm
$$ \norm{dP}{C^3_b} = \norm{dP}{C_b(N)} + \norm{d^2P}{C_b(N)} + \norm{d^3P}{C_b(N)} + \norm{d^4P}{C_b(N)}. $$
We now estimate, using \eqref{one} and \eqref{two},
\begin{align*}
\lVert d^2P_u( &(\nabla u)^2, (I-P_u) (\partial_t^2u + \Delta^2 u)) \Delta u_t\rVert_{L^1} \\ &\lesssim
\norm{\Delta u_t}{L^2} \norm{\nabla u}{L^{\infty}}^2 \bigg[\norm{u_t}{L^{\infty}} \norm{u_t}{L^2} + \norm{\nabla^2 u}{L^{\infty}}\norm{\Delta u}{L^2} + \norm{\nabla u}{L^{\infty}}^2 \norm{\Delta u}{L^2}\\
& + \norm{\nabla \Delta u}{L^2} \norm{\nabla u}{L^{\infty}} + \norm{\nabla u}{L^4}^2 \norm{\nabla u}{L^{\infty}}^2\bigg],
\end{align*}
\begin{align*}
\lVert dP_u( &\Delta u, (I-P_u) (\partial_t^2u + \Delta^2 u)) \Delta u_t\rVert_{L^1}\\
& \lesssim \norm{\Delta u_t}{L^2} \norm{\Delta u}{L^{\infty}} \bigg[\norm{u_t}{L^{\infty}} \norm{u_t}{L^2} + \norm{\nabla^2 u}{L^{\infty}}\norm{\Delta u}{L^2} + \norm{\nabla u}{L^{\infty}}^2 \norm{\Delta u}{L^2}\\
& + \norm{\nabla \Delta u}{L^2} \norm{\nabla u}{L^{\infty}} + \norm{\nabla u}{L^4}^2 \norm{\nabla u}{L^{\infty}}^2\bigg]\\
& = \norm{\Delta u_t}{L^2} \norm{\Delta u}{L^{\infty}} \big[\norm{u_t}{L^{\infty}} \norm{u_t}{L^2} + \norm{\nabla^2 u}{L^{\infty}}\norm{\Delta u}{L^2}\big]\\  
&~~~ + h(t)^2\norm{\Delta u_t}{L^2} \norm{\Delta u}{L^{\infty}}\big[ \norm{\Delta u}{L^2} + \norm{\nabla u}{L^{4}}^2\big]\\
&~~~ + h(t)\norm{\Delta u_t}{L^2} \norm{\Delta u}{L^{\infty}} \norm{\nabla \Delta u}{L^2},
\end{align*}
where we set $ h(t) := \norm{\nabla u(t)}{L^{\infty}}$. We further note that the equality is up to the constant from the estimate and hence proceed by estimating
\begin{align*}
\lVert dP_u(& \nabla  u, \nabla [(I-P_u) (\partial_t^2u + \Delta^2 u)]) \Delta u_t\rVert_{L^1} \\ 
&\lesssim \norm{\Delta u_t}{L^2} \norm{\nabla  u}{L^{\infty}} \bigg[\norm{u_t}{L^{\infty}} \norm{\nabla u_t}{L^2} + \norm{\nabla  u}{L^2} \norm{u_t}{L^{\infty}}^2 +\norm{\Delta^2 u}{L^2}\norm{\nabla u}{L^{\infty}}\\
& +  \norm{\nabla \Delta  u}{L^2}(\norm{\nabla^2 u}{L^{\infty}} + \norm{\nabla u}{L^{\infty}}^2) + \norm{\nabla u}{L^{\infty}} \norm{\Delta u}{L^2} \norm{\nabla^2 u}{L^{\infty}}+ \norm{ \Delta u}{L^2} \norm{\nabla u}{L^{\infty}}^3\\
& + \norm{\nabla u}{L^2} \norm{\nabla u}{L^{\infty}}^4\bigg].
\end{align*}
The latter upper bound equals the sum of
$$h(t)\norm{\Delta u_t}{L^2} \big[\norm{u_t}{L^{\infty}} \norm{\nabla u_t}{L^2} + \norm{\nabla  u}{L^2} \norm{u_t}{L^{\infty}}^2 +  \norm{\nabla \Delta  u}{L^2}\norm{\nabla^2 u}{L^{\infty}} + \norm{ \Delta u}{L^2} \norm{\nabla u}{L^{\infty}}^3\big],$$
and 
$$ h^2(t)\norm{\Delta u_t}{L^2} \big[ \norm{\Delta^2 u}{L^2} + \norm{\nabla u}{L^{\infty}} \norm{\nabla \Delta u}{L^2} + \norm{\Delta u}{L^2} \norm{\nabla^2 u}{L^{\infty}} + \norm{\nabla u}{L^2} \norm{\nabla u}{L^{\infty}}^3\big].$$
We calculate
\begin{align*}
&\nabla [ dP_u(\Delta u, u_t) + 2 dP_u(\nabla u, \nabla u_t) + d^2P_u((\nabla u)^2, u_t)]\\
&= d^2P_u(\nabla u, \Delta u, u_t) + dP_u(\nabla \Delta u, u_t) + dP_u(\Delta u, \nabla  u_t) + 2 d^2P_u((\nabla u)^2, \nabla u_t)\\
&~~~+2 dP_u(\nabla^2 u, \nabla u_t) + 2 dP_u(\nabla u, \nabla^2 u_t) + d^3P_u((\nabla u)^3, u_t)\\
&~~~ + 2d^2P_u(\nabla u, \nabla^2 u, u_t) + d^2P_u((\nabla u)^2, \nabla u_t),
\end{align*}
and hence
\begin{align*}
\lVert \nabla [& dP_u(\Delta u, u_t) + 2 dP_u(\nabla u, \nabla u_t) + d^2P_u((\nabla u)^2, u_t)]\cdot  \nabla [ (I-P_u) (\partial_t^2u + \Delta^2 u)]\rVert_{L^1}\\
& \lesssim \big(\norm{\Delta u}{L^2} \norm{\nabla  u}{L^{\infty}}\norm{u_t}{L^{\infty}} + \norm{\nabla \Delta u}{L^2} \norm{  u_t}{L^{\infty}} + (\norm{\Delta u}{L^{\infty}} + \norm{\nabla u}{L^{\infty}}^2) \norm{ \nabla  u_t}{L^2}\\
&~~~ + \norm{\Delta u_t}{L^2} \norm{  \nabla u}{L^{\infty}} + \norm{\nabla u}{L^{\infty}}^3 \norm{  u_t}{L^2}   \big) \bigg[\norm{u_t}{L^{\infty}} \norm{\nabla u_t}{L^2} + \norm{\nabla  u}{L^2} \norm{u_t}{L^{\infty}}^2\\
&~~~+\norm{\Delta^2 u}{L^2}\norm{\nabla u}{L^{\infty}} +  \norm{\nabla \Delta  u}{L^2}(\norm{\nabla^2 u}{L^{\infty}} + \norm{\nabla u}{L^{\infty}}^2) + \norm{\nabla u}{L^{\infty}} \norm{\Delta u}{L^2} \norm{\nabla^2 u}{L^{\infty}}\\
&~~~+ \norm{ \Delta u}{L^2} \norm{\nabla u}{L^{\infty}}^3 + \norm{\nabla u}{L^2} \norm{\nabla u}{L^{\infty}}^4\bigg].
\end{align*}
We now collect all terms which are quadratic, linear or constant in $h(t)$, i.e. the latter bound equals
$$ J_1(u) + h(t)J_2(u) + h(t) J_3(u) + h^2(t) J_4(u),$$
where
\begin{align*}
J_1(u) &=(\norm{\nabla \Delta u}{L^2}\norm{u_t}{L^{\infty}}+\norm{\Delta u}{L^{\infty}}\norm{\nabla u_t}{L^2})\big[\norm{u_t}{L^{\infty}} \norm{\nabla u_t}{L^2} + \norm{\nabla  u}{L^2} \norm{u_t}{L^{\infty}}^2\\
&~~~ +  \norm{\nabla \Delta  u}{L^2}\norm{\nabla^2 u}{L^{\infty}} + \norm{ \Delta u}{L^2} \norm{\nabla u}{L^{\infty}}^3\big],\\[3pt]
J_2(u) &= (\norm{\nabla \Delta u}{L^2}\norm{u_t}{L^{\infty}}+\norm{\Delta u}{L^{\infty}}\norm{\nabla u_t}{L^2}) \big[\norm{\Delta^2 u}{L^2} + \norm{\nabla u}{L^{\infty}} \norm{\nabla \Delta u}{L^2}\\
&~~~ + \norm{\Delta u}{L^2} \norm{\nabla^2 u}{L^{\infty}} + \norm{\nabla u}{L^2} \norm{\nabla u}{L^{\infty}}^3\big],\\[3pt]
J_3(u) &=  (\norm{\Delta u}{L^2}\norm{u_t}{L^{\infty}} + \norm{\nabla u_t}{L^2} \norm{\nabla u}{L^{\infty}}+ \norm{\Delta u_t}{L^2} + \norm{\nabla u}{L^{\infty}}^2 \norm{u_t}{L^2}) \big[\norm{u_t}{L^{\infty}} \norm{\nabla u_t}{L^2}\\[3pt]
&~~~ + \norm{\nabla  u}{L^2} \norm{u_t}{L^{\infty}}^2 +  \norm{\nabla \Delta  u}{L^2}\norm{\nabla^2 u}{L^{\infty}}+ \norm{ \Delta u}{L^2} \norm{\nabla u}{L^{\infty}}^3\big],\\[3pt]
J_4(u) &= (\norm{\Delta u}{L^2}\norm{u_t}{L^{\infty}} + \norm{\nabla u_t}{L^2} \norm{\nabla u}{L^{\infty}}+ \norm{\Delta u_t}{L^2} + \norm{\nabla u}{L^{\infty}}^2 \norm{u_t}{L^2}) \big[ \norm{\Delta^2 u}{L^2}\\
&~~~ + \norm{\nabla u}{L^{\infty}} \norm{\nabla \Delta u}{L^2} + \norm{\Delta u}{L^2} \norm{\nabla^2 u}{L^{\infty}} + \norm{\nabla u}{L^2} \norm{\nabla u}{L^{\infty}}^3\big].
\end{align*}
We note that the energy is conserved, ie. for $ t \in [0, T) $
\begin{equation}\label{energy-conservation}
2E(u(t)) = \norm{\Delta u(t)}{L^2}^2 + \norm{\partial_t u(t)}{L^2}^2 = \norm{\Delta u_0}{L^2}^2 + \norm{u_1}{L^2}^2 = 2E(u_0, u_1),
\end{equation}
and further, this implies the bounds
\begin{align}
&\sup_{t \in [0, T)}\norm{\nabla u(t)}{L^2} \lesssim \sqrt{1 + T}(\sqrt{E(u_0, u_1)} + \norm{\nabla u_0}{L^2}),~~\text{and}\label{this1}\\
& \sup_{t \in [0, T)}\norm{ u(t)- u_0}{L^2} \lesssim T\sqrt{E(u_0, u_1)}.\label{this2}
\end{align}
We  recall the following cases of Gagliardo-Nirenberg's interpolation for $ n = 2$
\begin{align}
\norm{\Delta u}{L^{\infty}}+ \norm{\nabla \Delta u}{L^2} &\lesssim \norm{\Delta^2u}{L^2}^{\f{1}{2}}\norm{\Delta u}{L^2}^{\f{1}{2}},~~\norm{u_t}{L^{\infty}} \lesssim \norm{\Delta u_t}{L^2}^{\f{1}{2}}\norm{ u_t}{L^2}^{\f{1}{2}},~~ \label{eins}\\
\norm{\nabla u}{L^{\infty}} &\lesssim \norm{\Delta^2 u}{L^2}^{\f{1}{3}}\norm{ \nabla u}{L^2}^{\f{2}{3}},~~\norm{\nabla u}{L^{4}} \lesssim \norm{\Delta^2 u}{L^2}^{\f{1}{6}}\norm{ \nabla u}{L^2}^{\f{5}{6}},~~\text{and}\label{zwei}\\
\norm{\nabla u_t}{L^{4}} &\lesssim \norm{\Delta u_t}{L^2}^{\f{3}{4}}\norm{ u_t}{L^2}^{\f{1}{4}}.\label{drei}
\end{align}
Setting 
$$ \mathcal{E}(u(t)) := \norm{\Delta u_t(t)}{L^2} + \norm{\Delta^2 u(t)}{L^2},~ t \in [0, T),$$
by \eqref{eins}, \eqref{zwei}  and the estimates above, there exists a constant $C(T) = C(N, u_0, u_1) ( 1+ T)^{\alpha}$ for some $ \alpha > 0 $, such that $ C(N, u_0, u_1)$ only depends on the norm $ \norm{dP}{C^3_b} $, the optimal Sobolev constant in Gagliardo-Nirenberg's interpolation and $ E(u_0, u_1),~ \norm{\nabla u_0}{L^2}$ such that the following holds.
\begin{align}
\f{d}{dt} \mathcal{E}^2(u(t)) &\leq C(T) ( 1 + h(t) + h^2(t) )(\mathcal{E}(t) +\mathcal{E}^2(t))\\\nonumber
& \leq C(T)( 1 + h^2(t) )(1 +\mathcal{E}^2(t)),~ t \in [0, T).
\end{align}
Using the idea from \citep{fan2010regularity}, we now apply the sharp Sobolev inequality of Brezis-Gallouet-Wainger from \citep{brezis1980note} in order to bound (we assume $ u$ is not a constant)
\begin{equation}
 h(t) \leq \tilde{C}\norm{\nabla u(t)}{H^1}\left( 1 + \log^{\f{1}{2}}\left(1 + \f{\norm{\nabla u(t)}{H^2}^2}{\norm{\nabla u(t)}{H^1}^2}\right) \right),~~ t \in [0, T).
\end{equation}
Thus, using \eqref{eins},~\eqref{this1} and \eqref{energy-conservation}, 
\begin{equation}
h^2(t) \leq C(T)\left( 1 + \log\left(1 + \mathcal{E}^2(t)\right) \right),~~ t \in [0, T),
\end{equation}
and hence
\begin{align}\label{Gron}
\f{d}{dt}( e +  \mathcal{E}^2(u(t))) &\leq C(T)  \log\left(e + \mathcal{E}^2(t)\right)(e  + \mathcal{E}^2(t)),~ t \in [0, T).
\end{align}
 This suffices for a Gronwall-type inequality for $ \log( e + \mathcal{E}^2(t))$ and hence by \eqref{energy-conservation} and \eqref{eins}, \eqref{zwei} and \eqref{drei}, we have
 $$ \limsup_{t \to T} ( \norm{u_t}{H^2}^2 + \norm{\nabla u}{H^3}^2) < \infty,$$
 as long as $ T < \infty$.\\[5pt]
\textit{\underline{Case: n = 1}}~~Here, by Gagliardo-Nirenberg's estimate, we infer the bound
\begin{align}\label{easy}
\norm{\nabla u}{L^{\infty}} \lesssim \norm{\nabla^2 u}{L^2}^{\f{1}{2}}\norm{\nabla u}{L^2}^{\f{1}{2}}.
\end{align} 
Hence, the a priori bound is derived similarly  for $ (\nabla u(t), u_t(t)) \in H^2(\R) \times H^1(\R)$. We note  
\begin{align*}
\f{d}{dt} \int_{\R} &| \nabla u_t|^2 + |\nabla \Delta u|^2~dx = -\int_{\R} dP_u( \nabla u, (I-P_u)( \partial_t^2 u + \Delta^2 u)) \cdot\nabla u_t\\
 &~~- \int_{\R} (I-P_u)( \partial_t^2 u + \Delta^2 u)) \cdot (d^2P_u((\nabla u)^2,  u_t) + dP_u(\nabla^2 u, u_t) + dP_u(\nabla u, \nabla u_t))dx.
\end{align*}
Thus we estimate, as before
\begin{align*}
\lVert dP_u( \nabla u, (I-P_u)& (\partial_t^2u + \Delta^2 u)) \nabla u_t\rVert_{L^1}\\
& \lesssim \norm{\nabla u_t}{L^2} \norm{\nabla u}{L^{\infty}} \bigg[\norm{u_t}{L^{\infty}} \norm{u_t}{L^2} + \norm{\nabla^2 u}{L^{\infty}}\norm{\nabla^2 u}{L^2} + \norm{\nabla u}{L^{\infty}}^2 \norm{\nabla^2 u}{L^2}\\
& + \norm{\nabla^3 u}{L^2} \norm{\nabla u}{L^{\infty}} + \norm{\nabla u}{L^4}^2 \norm{\nabla u}{L^{\infty}}^2\bigg],~~\text{and}
\end{align*}
\begin{align*}
\lVert &  (d^2P_u((\nabla u)^2,  u_t) + dP_u(\nabla^2 u, u_t) + dP_u(\nabla u, \nabla u_t)) [(I-P_u) (\partial_t^2u + \Delta^2 u)]\rVert_{L^1} \\ 
&\lesssim (\norm{ u_t}{L^2} \norm{\nabla  u}{L^{\infty}}^2 + \norm{\nabla^2 u}{L^2} \norm{u_t}{L^{\infty}} + \norm{\nabla u_t}{L^2}\norm{\nabla u}{L^{\infty}})
\bigg[\norm{u_t}{L^{\infty}} \norm{u_t}{L^2}\\
&~~~ + \norm{\nabla^2 u}{L^{\infty}}\norm{\nabla^2 u}{L^2} + \norm{\nabla u}{L^{\infty}}^2 \norm{\nabla^2 u}{L^2} + \norm{\nabla^3 u}{L^2} \norm{\nabla u}{L^{\infty}} + \norm{\nabla u}{L^4}^2 \norm{\nabla u}{L^{\infty}}^2\bigg]. 
\end{align*}
Hence from the interpolation estimates \eqref{easy},
\begin{align}
\norm{\nabla^2u }{L^{\infty}} \lesssim \norm{\nabla^4u}{L^2}^{\f{1}{4}}\norm{\nabla^2u}{L^2}^{\f{3}{4}},~~\norm{u_t}{L^{\infty}} \lesssim \norm{\nabla u_t}{L^2}^{\f{1}{2}}\norm{u_t}{L^2}^{\f{1}{2}},
\end{align}
and \eqref{energy-conservation},~\eqref{this1}, there holds (for $C(T) > 0$ as before) 
\begin{equation}
\f{d}{dt}(1 + \mathcal{E}(t)) \leq C(T)( 1 + \mathcal{E}(t)),~ t \in [0, T)
\end{equation}
which suffices to use a Gronwall argument in order to conclude the proof.
\section{A uniqueness argument}\label{uniqueness}
We now give a short argument for the uniqueness of solutions $ u : \R^n \times [0, T) \to N,~ n = 1,2,3 $ with
 \begin{equation}\label{class}
 u- u(0) \in C^0([0,T), H^4(\R^n))\cap C^1([0,T), H^2(\R^n)).
 \end{equation}
 Setting $ w = u-v $ for solutions $ u, ~v $ of \eqref{CP} in the class \eqref{class} with $ u(0) = v(0),~ u_t(0) = v_t(0)$, we provide a Gronwall type argument in the energy space, i.e. more precisely for the norm $ \norm{w_t}{L^2}^2 + \norm{ w}{H^2}^2 $. We note the interpolation estimate 
 $$ \norm{w}{L^{\infty}} \lesssim \norm{\Delta w}{L^2}^{\f{n}{4}} \norm{w}{L^2}^{\f{4-n}{4}},~~ n = 1,2,3,$$
 and the identity
 \begin{align}\label{energy}
 \f{d}{2 dt}\left(\int_{\R^n}|w_t|^2 + |\Delta w|^2\right)~dx = I_1 + I_2 + I_3,
 \end{align}
 where
 \begin{align*}
 &I_1 = \int_{\R^n} \big[ dP_u( u_t u_t + 4 \nabla u \cdot \nabla \Delta u + \Delta u \Delta u + 2 \nabla^2 u \cdot \nabla^2 u) + d^3P_u( \nabla u)^4\\
 &~~~~~+ d^2P_u(2 (\nabla u)^2 \Delta u + 4(\nabla u)^2 \cdot \nabla^2 u ) \big]( P_v - P_u) w_t~dx\\
 &~~~~~+ \int_{\R^n} \big[ (dP_u - dP_v)( u_t u_t + 4 \nabla u \cdot \nabla \Delta u + \Delta u \Delta u + 2 \nabla^2 u \cdot \nabla^2 u)\\
  &~~~~~+ (d^3P_u - d^3P_v)( \nabla u)^4 + (d^2P_u - d^2P_v)(2 (\nabla u)^2 \Delta u + 4(\nabla u)^2 \cdot \nabla^2 u ) \big](P_u - P_v) u_t~dx\\
  &I_2 = \int_{\R^n}\big[ dP_v( u_t w_t + w_t v_t + 4 \nabla w \cdot \nabla \Delta u + \Delta w \Delta u + \Delta v \Delta w\\
  &~~~~~ + 2 \nabla^2 w \cdot \nabla^2 u + 2 \nabla^2 v \cdot \nabla^2 w) + d^3P_v(\nabla w \cdot (\nabla u)^3 + \nabla w \cdot (\nabla u)^2 \nabla v +  \nabla w \cdot (\nabla v)^2 \nabla u +  \nabla w \cdot (\nabla v)^3  )\\
  &~~~~~ + d^2P_v( 2 \nabla w \cdot \nabla u \Delta u + 2 \nabla v \cdot \nabla w \Delta u + 2 (\nabla v)^2  \Delta w +  4 \nabla w \cdot \nabla u \cdot \nabla^2 u\\
  &~~~~~ + 4 \nabla v \cdot \nabla w \cdot  \nabla^2 u + 4 (\nabla v)^2 \cdot \nabla^2 w )\big] (P_u - P_v) u_t~dx\\
  & I_3 = \int_{\R^n } 4dP_v(\nabla v, \nabla \Delta w)(P_u - P_v) u_t~dx.
 \end{align*}
This follows from 
 \begin{align*}
  (I-P_u )&(\partial_t^2 u  + \Delta^2 u) - (I-P_v )(\partial_t^2 v + \Delta^2 v)\\
  &= (P_v - P_u)[ (I-P_u)(\partial_t^2 u + \Delta^2 u)]\\
  &~~~ +  (I-P_v )[(I-P_u)(\partial_t^2 u + \Delta^2 u) - (I-P_v )(\partial_t^2 v + \Delta^2 v)],
 \end{align*}
 and 
 $$ (I-P_v)w_t = (I-P_v)u_t = (P_u - P_v)u_t.$$
 We further note 
 \begin{align*}
 \int_{\R^n}dP_v(& \nabla v, \nabla \Delta w) (P_u - P_v) u_t~dx = - \int_{\R^n} \left[d^2P_v( (\nabla v)^2,  \Delta w) + dP_v(\Delta v, \Delta w)\right] (P_u - P_v) u_t~dx\\
 &~~~ - \int_{\R^n}dP_v( \nabla v, \Delta w) (P_u - P_v) \nabla u_t + dP_v (\nabla v, \Delta w) (dP_u - dP_v) (\nabla u, u_t)~dx\\
 &~~~ - \int_{\R^n} dP_v( \nabla v, \Delta w ) dP_v (\nabla w, u_t)~dx.
 \end{align*}
 Hence, we estimate
 \begin{align*}
 I_1 &\lesssim (\norm{w_t}{L^2}\norm{w}{L^{\infty}} + \norm{w}{L^{\infty}}^2\norm{u_t}{L^2})( \norm{u_t}{L^2}\norm{u_t}{L^{\infty}} + \norm{\nabla \Delta u}{L^2}\norm{\nabla u}{L^{\infty}} + \norm{\nabla u}{L^2}\norm{\nabla u}{L^{\infty}}^3\\
 &~~~~~ + \norm{\Delta u}{L^2} \norm{\nabla u}{L^{\infty}}^2 + \norm{\Delta u}{L^2} \norm{ \nabla^2 u}{L^{\infty}}),
 \end{align*}
 \begin{align*}
 I_2 &\lesssim \norm{w}{L^{\infty}}\norm{u_t}{L^2}( \norm{\nabla w}{L^2} \norm{\nabla u }{L^{\infty}}^3 + \norm{\nabla w}{L^2} \norm{\nabla v }{L^{\infty}} \norm{\nabla u }{L^{\infty}}^2 + \norm{\nabla w}{L^2} \norm{\nabla v }{L^{\infty}}^2 \norm{\nabla u }{L^{\infty}}\\[4pt]
 &~~~ + \norm{\nabla w}{L^2} \norm{\nabla v }{L^{\infty}}^3 + \norm{\nabla w}{L^2}\norm{\nabla u}{L^{\infty}}\norm{\nabla^2 u}{L^{\infty}} +  \norm{\nabla w}{L^2}\norm{\nabla v}{L^{\infty}}\norm{\nabla^2 u}{L^{\infty}}\\[4pt]
 &~~~ + \norm{\nabla v}{L^{\infty}}^2\norm{\Delta w}{L^2} + \max\{ \norm{u_t}{L^{\infty}}, \norm{v_t}{L^{\infty}}\}\norm{w_t}{L^2} + \max\{ \norm{\nabla^2 u}{L^{\infty}}, \norm{\nabla^2 v}{L^{\infty}}\}\norm{\Delta w}{L^2}) \\
 &~~~+\norm{w}{L^{\infty}}\norm{u_t}{L^{\infty}}\norm{\nabla \Delta u}{L^2}\norm{\nabla w}{L^2} 
 \end{align*}
 and
 \begin{align*}
 &I_3 \lesssim  \norm{w}{L^{\infty}}\norm{\nabla u_t}{L^2}\norm{\Delta w}{L^2}\norm{\nabla v}{L^{\infty}} + \norm{w}{L^{\infty}}\norm{u_t}{L^2} \norm{\Delta w}{L^2}\norm{\nabla v}{L^{\infty}}\norm{\nabla u}{L^{\infty}}\\[4pt]
 &~~~~~~+ \norm{\Delta w}{L^2} \norm{\nabla w}{L^2}\norm{u_t}{L^{\infty}} \norm{\nabla v}{L^{\infty}} + \norm{w}{L^{\infty}}\norm{u_t}{L^2}\norm{\Delta w}{L^2} \norm{\nabla v}{L^{\infty}}^2\\
 &~~~~~~+ \norm{w}{L^{\infty}}\norm{u_t}{L^{\infty}}\norm{\Delta w}{L^2} \norm{\Delta v }{L^2}.
 \end{align*}
 We set 
 $$ \mathcal{E}^2(t) : = \norm{w_t}{L^2}^2 + \norm{ w}{H^2}^2.$$
 ~~\\
 Using the aforementioned interpolation inequality, we obtain in particular $ \norm{w}{L^{\infty}} \lesssim \mathcal{E}(t)$. Since also
 \begin{align}
& \f{d}{dt}\int_{\R^n}|\nabla w|^2 ~dx \leq  \norm{w_t}{L^2}^2 + \norm{\Delta w}{L^2}^2 \leq \mathcal{E}^2(t),~~\text{and}\\
&\f{d}{dt}\int_{R^n} |w|^2 ~dx \leq \norm{w}{L^2}^2 + \norm{w_t}{L^2}^2 \leq \mathcal{E}^2(t),
 \end{align}
 estimating \eqref{energy} gives
 \begin{align}
 \f{d}{dt}\mathcal{E}(t)^2 \lesssim ( 1 + \norm{\nabla u}{H^{3}}^4 + \norm{u_t}{H^{2}}^4 + \norm{\nabla v}{H^{3}}^4 + \norm{v_t}{H^{2}}^4)\mathcal{E}^2(t) =: C(u,v) \mathcal{E}^2(t).
 \end{align}
This suffices for uniqueness, as long as $ C(u,v)$ stays bounded in time.
 We also remark that in $ n = 1 $, in order to conclude uniqueness from similar arguments, it suffices for a smooth solution $u$ to stay bounded in $ u(t) \in H^3(\R),~ \partial_t u(t) \in H^1(\R)$.

 \bibliography{mybib}

\begin{thebibliography}{10}

\bibitem{brezis1980note}
Ha{\"\i}m Br{\'e}zis and Stephen Wainger.
\newblock A note on limiting cases of sobolev embeddings and convolution
  inequalities.
\newblock {\em Communications in Partial Differential Equations},
  5(7):773--789, 1980.

\bibitem{fan2010regularity}
Jishan Fan and Tohru Ozawa.
\newblock On regularity criterion for the 2{D} wave maps and the 4{D}
  biharmonic wave maps.
\newblock In {\em Current advances in nonlinear analysis and related topics},
  volume~32 of {\em GAKUTO Internat. Ser. Math. Sci. Appl.}, pages 69--83.
  Gakkotosho, Tokyo, 2010.

\bibitem{grillakisgeba}
Dan-Andrei Geba and Manoussos~G Grillakis.
\newblock {\em An introduction to the theory of wave maps and related geometric
  problems}.
\newblock World Scientific, 2017.

\bibitem{jiang}
Jiang Guoying.
\newblock 2-harmonic maps and their first and second variational formulas.
\newblock {\em Note di Matematica}, 28(supn1):209--232, 2008.

\bibitem{LammSchnaubeltHerrSchmid}
Sebastian Herr, Tobias Lamm, Tobias Schmid, and Roland Schnaubelt.
\newblock Biharmonic wave maps: local wellposedness in high regularity.
\newblock {\em Nonlinearity}, 33(5):2270, 2020.

\bibitem{LammSchnaubeltHerr}
Sebastian Herr, Tobias Lamm, and Roland Schnaubelt.
\newblock Biharmonic wave maps into spheres.
\newblock {\em Proceedings of the American Mathematical Society},
  148(2):787--796, 2020.

\bibitem{schmid2}
Tobias Schmid.
\newblock Global results for a cauchy problem related to biharmonic wave maps.
\newblock {\em arXiv preprint, arXiv:2102.12881}, 2021.

\bibitem{Sschmid}
Tobias Schmid.
\newblock Local wellposedness and global regularity results for biharmonic wave
  maps.
\newblock \url{http://dx.doi.org/10.5445/IR/1000128147}, (2021).

\bibitem{shatah1988weak}
Jalal Shatah.
\newblock Weak solutions and development of singularities of the su (2)
  $\sigma$-model.
\newblock {\em Communications on pure and applied mathematics}, 41(4):459--469,
  1988.

\bibitem{shatahstruwe}
Jalal Shatah and Michael Struwe.
\newblock {\em Geometric wave equations}.
\newblock American Mathematical Society, Providence, RI, 1998.

\bibitem{Strzelecki}
Pawel Strzelecki.
\newblock On biharmonic maps and their generalizations.
\newblock {\em Calculus of Variations and Partial Differential Equations},
  18(4):401--432, Dec 2003.

\end{thebibliography}

\end{document}